\documentclass[12pt]{article}
       \usepackage{amsmath,amsthm}
       \usepackage{amssymb,amsbsy}
       \usepackage{amsfonts}

\usepackage{cite}
\usepackage{times}
\newcommand{\footnotep}[1]{\protect\footnote{#1}}

\begin{document}

\author{
{\large {\bf  \framebox{W. Marcinek}
}}\\[3pt] {\it Institute of
Theoretical Physics, University of Wroc\l aw}\\ Pl. Maxa Borna 9,
50-204 Wroc{\l}aw, Poland
}

\title{{\bf PARTICLES OF GENERALIZED STATISTICS, QUANTUM LOGIC AND CATEGORIES}
\thanks{These are last unfinished ideas by W\l adys\l aw Marcinek
who unexpectedly passed away on June 9, 2003.
Prepared for publication by Steven Duplij.}
\thanks{Published in the W. Marcinek's Memorial volume of the
Kharkov University Journal Vestnik KNU, ser. "Nuclei, particles and fields".
- 2003. - V. 601. - N 2(22). - p. 135--146.}
 }
\thispagestyle{empty}
 \date{May 25, 2003}
\maketitle
\begin{abstract}
A collection of brilliant and original unfinished ideas by Wladyslaw
Marcinek (1952-2003) in particle interactions, categorical approach to
generalized statistics, qubits and quantum logic, entwined operators, cobordisms and
noncommutative Fock space.
\end{abstract}

\newpage
\smallskip
\section{Composite Systems and Regularity
\footnotep{Last change May 19, 2003}}

Let us assume that an interacting system of charged particles and quanta is
given as a starting point for the study of particle processes. It is natural
to assume that there is a set of initial configurations of the system and
there is also a set of final configurations representing possible results of
interactions. An arbitrary initial configuration can be transform into final
one as a result of a sequences of vertex interactions. Every such sequence is
said to be a \textit{process}. The proper physical meaning have the initial
and final configurations, intermediate steps are virtual. A process with an
initial configuration and one final configuration and without intermediate
steps is said to be \textit{primary}. An arbitrary process can be built from
primary ones. If the time evolution of the system can be described in an
unique way as a sequence of transformations of initial configuration, then we
say that the system is \textit{consistent} and equipped with a
\textit{generalized quantum history}. In this approach to quantum theory the
whole universe is represented by a class of 'histories'. In this formalism the
standard Hamiltonian time-evolution is replaced by a partial semigroup called
a 'temporal support' \cite{gel-hart}.

It is natural to assume that the whole world is divided into two parts: a
classical particle system and its quantum environment. The classical system
represents an observed reality, this particles which really exist and can be
detect. The quantum environment represents all quantum possibilities to become
a part of the reality in the future \cite{haa}.

An algebraic model of an interacting system of charged particles and quanta of
an external field has been developed by the author \cite{wmq,wmb,gqsi}. The
model is based on the assumption that a given charged particle is transform
under interaction into a composite system consisting a charge and quanta of
the quantum field. The system containing a particle dressed with a single
quantum of external field is said to be a \textit{quasiparticle}.
Quasiparticles behave like free particles equipped with generalized
statistics. In this paper we develop the algebraic model considered in
previous publications for systems with generalized statistics
\cite{wmq,castat,gqs}. We use the category theory in order to construct the
model. The initial data for construction of category $\mathfrak{C}$ relevant
for our model is described. In our approach objects of the category
$\mathfrak{C}$ represent physical objects like charged particles,
quasiparticles, quasiholes, different species of quanta of an external field,
etc... Functors or multifunctors represent interactions, creations and
annihilations, etc... There is a collection of (multi)-functors representing
primary act of interactions. Starting from the collection we can construct a
collection of higher - ary functors representing multiquasiparticle processes.
The problem is to describe a collection of natural isomorphisms such that our
constructions are unique up to these isomorphisms.

Let us consider a system of charged particles interacting with a quantum
environment. States of the system represent all initial configurations. These
initial configurations are transformed into final configuration as a result of
interaction of the system with the quantum environment. In this transformation
every charged particle is converted into a composite system consisting a few
charges and quanta from the quantum environment. Such a transformation is said
to be a composition. A conversion of a composite system into his parts is
called a decomposition. Note that both transformations, a composition and a
decomposition need not to be invertible but they can be regular.

A system which contains a charge and certain number of quanta as a result of
interaction with the quantum field is said to be a \textit{dressed particle}
\cite{gqsi}. Next we assume that every dressed particle is a composite object
equipped with an internal structure. Obviously the structure of dressed
particles is determined by the interaction with the quantum field. We describe
a dressed particle as a non-local system which contains $n$ centers
(vertexes). Two systems with $n$ and $m$ centers, respectively, can be
'composed' into one system with $n+m$ centers. All centers as members of a
given system behave like free particles moving on certain effective space
$\mathcal{M}$. Every center is also equipped with ability for absorption and
emission of quanta of the quantum field. Assume that there are quanta of $N$
different kindS. A centre dressed with a single quantum of the field is said
to be a \textit{quasiparticle}. In this way $N$ is not a number of real
particles but a number of quasparticles. In our approach a center equipped
with two quanta forms a system of two quasiparticles.

A center with an empty place for a single quantum is said to be a
\textit{quasihole}. A centre which contains any quantum is said to be
\textit{neutral}. A neutral center can be transform into a quasiparticle or a
quasihole by an absorption or emission process of single quantum,
respectively. In this way the process of absorption of quanta of quantum field
by a charged particle is equivalent to a creation of quasiparticles and
emission provide to annihilation of quasiparticles. Note that there is also
the process of mutual annihilation of quasiparticles and quasiholes.

\section{Generalized Qubits and Quantum Statistics
\footnotep{Last change May 7, 2003}}

In the last years a few different approaches to quantum statistics which
generalize the usual boson or fermion statistics has been intensively
developed by several authors. The so--called $q$--statistics and corresponding
$q$--relations have been studied by Greenberg \cite{owg,gre}, Mohapatra
\cite{moh}, Fivel \cite{fi} and many others, see \cite{zag,mepe,bs2} for
example. The deformation of commutation relations for bosons and fermions
corresponding to quantum groups $SU_{q}(2)$ has been given by Pusz and
Woronowicz \cite{P,PW}. The $q$--relations corresponding to superparticles has
been considered by Chaichian, Kulisch and Lukierski \cite{ckl}. Quantum
deformations have been also studied by Vokos \cite{V}, Fairle and Zachos
\cite{FZ} and many others.

It is interesting that in the last years new and highly organized structures
of matter has been discovered. For example in fractional quantum Hall effect a
system with well defined internal order has appears \cite{zee}. Another
interesting structures appear in the so called $\frac{1}{2}$ electronic
magnetotransport anomaly \cite{jai,dst}, high temperature superconductors or
laser excitations of electrons. In these cases certain anomalous behaviour of
electron have appear. The existence of new ordered structures depends on the
existence of certain specific additional excitations. Hence we can restrict
our attention to the study of possibility of appearance for these excitations.
For the description of such possible excited states we use the concept of
dressed particles. We assume that every charged particle is equipped with
ability to absorb quanta of the external field. A system which contains a
particle and certain number of quanta as a result of interaction with the
external field is said to be dressed particle. A particle without quantum is
called undressed or a quasihole. The particle dressed with two quanta of
certain species is understand as a system of two new objects called
quasiparticles. A quasiparticle is in fact the charged particle dressed with a
single quantum. Two quasiparticles are said to be identical if they are
dressed with quanta of the same species. In the opposite case when the
particle is equipped with two different species of quanta then we have
different quasiparticles. We describe excited states as composition of
quasiparticles and quasiholes. It is interesting that quasiparticles and
quasiholes have also their own statistics. We give the following assumption
for the algebraic description of excitation spectrum of single dressed particle.

\textbf{Fundamental Assumptions:}

\emph{Assumption 0: The ground state.} There is a state $|0>=\mathbf{1}$
called the ground one. There is also the conjugate ground state $<0|\equiv
\mathbf{1}^{\ast}$. This is the state of the system before intersection.

\emph{Assumption 1: Elementary states.} There is an ordered (finite) set of
single quasiparticle states $S:=\{x^{i}:i=1,\ldots,N<\infty\}$. These states
are said to be elementary (simple). They represent elementary excitations of
the system. We assume that the set $S$ of elementary states forms a basis for
a finite linear space $E$ over a field of complex numbers $\mathbb{C}$.

\emph{Assumption 2: Elementary conjugated states.} There is also a
corresponding set of single quasihole states $S^{\ast}:=\{x^{\ast
i}:i=N,N-1,\ldots,1\}$. These states are said to be conjugated. The set
$S^{\ast}$ of conjugate states forms a basis for the complex conjugate space
$E^{\ast}$. The pairing $(.|.):E^{\ast}\otimes E\longrightarrow\mathbb{C}$ is
given by $(x^{\ast i}|x^{j}):=\delta^{ij}$.

\emph{Assumption 3: Composite states.} There is a set of projectors $\Pi
_{n}:E^{\otimes n}\longrightarrow E^{\otimes n}$ such that we have a
$n$-multilinear mapping $\odot_{n}:E^{\times n}\longrightarrow E^{\otimes n}$
defined by the following formula $x^{i_{1}}\odot\cdots\odot x^{i_{n}}:=\Pi
_{n}(x^{i_{1}}\otimes\cdots\otimes x^{i_{n}})$. The set of $n$%
--multiquasiparticle states is denoted by $P^{n}(S)$. All such states are
result of composition (or clustering) of elementary ones. These states are
also called composite states of order $n$. They represent additional
excitations charged particle under interaction. In this way for
multiquasiparticle states we have the following set of states $P^{n}%
(S):=\{x^{\sigma}\equiv x^{i_{1}}\odot\cdots\odot x^{i_{n}}:\sigma
=(i_{1},\ldots,i_{n})\in I\}$. Here $I$ is a set of sequences of indices such
that the above set of states forms a basis for a linear space $\mathcal{A}%
^{n}$. We have $\mathcal{A}^{n}=Im(\Pi_{n})$. Obviously we have where
$\mathcal{A}^{0}\equiv\mathbf{1}\mathbb{C}$, $\mathcal{A}^{1}\equiv E$ and
$\mathcal{A}^{n}\subset E^{\otimes n}$.

\emph{Assumption 4: Composite conjugated states.} We also have a set of
projectors $\Pi_{n}^{\ast}:E^{\ast\otimes n}\longrightarrow E^{\ast\otimes n}$
and the corresponding set of composite conjugated states of length $n$
\begin{equation}
P^{n}(S^{\ast}):=\{x^{\ast\sigma}\equiv x^{\ast i_{n}}\odot\cdots\odot x^{\ast
i_{1}}:\sigma=(i_{1},\ldots,i_{n})\in I\}.
\end{equation}

The set $P^{n}(S^{\ast})$ of composite conjugated states of length $n$ forms a
basis for a linear space $\mathcal{A}^{\ast n}$.

\emph{Assumption 5: Algebra of states.} The set of all composite states of
arbitrary length is denoted by $P(S)$. For this set of states we have the
following linear space $\mathcal{A}:=\bigoplus\limits_{n}\ \mathcal{A}^{n}$.
If the formula
\begin{equation}
m(s\otimes t)s\equiv s\odot t:=\Pi_{m+n}(\tilde{s}\otimes\tilde{t})
\label{mul1}%
\end{equation}
for $s=\Pi_{m}(\tilde{s}),t=\Pi_{n}(\tilde{t})$, $\tilde{s}\in E^{\otimes n}$,
$\tilde{t}\in E^{\otimes m}$, defines an associative multiplication in
$\mathcal{A}$, then we say that we have an algebra of states. This algebra
represents excitation spectrum for single dressed particle.

\emph{Assumption 6: Algebra of conjugated states.} The set of composite
conjugated states of arbitrary length is denoted by $P(S^{\ast})$. We have
here a linear space $\mathcal{A}^{\ast}:=\bigoplus\limits_{n}\ \mathcal{A}%
^{\ast n}$. If $m$ is the multiplication in $\mathcal{A}$, then the
multiplication in $\mathcal{A}^{\ast}$ corresponds to the opposite
multiplication $m^{op}(t^{\ast}\otimes s^{\ast})=(m(s\otimes t))^{\ast}$.

We define creation operators for our model as multiplication in the algebra
$\mathcal{A}$ as $a_{s}^{+}t:=s\odot t,$ for $s,t\in\mathcal{A}$, where the
multiplication is given by (\ref{mul1}). For the ground state and annihilation
operators we assume that $\langle0|0\rangle=0$, $a_{s^{\ast}}|0\rangle=0$, for
$s^{\ast}\in\mathcal{A}^{\ast}.$The proper definition of action of
annihilation operators on the whole algebra $\mathcal{A}$ is a problem. For
the pairing $\langle-|-\rangle^{n}:\mathcal{A}^{\ast n}\otimes\mathcal{A}%
^{n}\longrightarrow\mathbb{C}$ we assume in addition that we have the
following formulae
\begin{equation}
\langle0|0\rangle^{0}:=0,\quad\langle i|j\rangle^{1}:=(x^{\ast i}%
|x^{j})=\delta^{ij},\ \ \ \langle s|t\rangle^{n}:=\langle\tilde{s}|P_{n}%
\tilde{t}\rangle n_{0}\quad\mbox{for}\quad n\neq2
\end{equation}
where $\tilde{s},,\tilde{t},\in E^{\otimes n}$, $P_{n}:E^{\otimes
n}\longrightarrow E^{\otimes n}$ is an additional linear operator and $$\langle
i_{1}\cdots i_{n}|j_{1}\cdots j_{n}\rangle_{0}^{n}:=\langle i_{1}|j_{1}%
\rangle^{1}\cdots\langle j_{n}|j_{n}\rangle^{1}.$$ Observe that we need two
sets $\Pi:=\{\Pi_{n}\}$ and $P:=\{P_{n}\}$ of operators and the action
\begin{equation}%
\begin{array}
[c]{l}%
a:s^{\ast}\otimes t\in\mathcal{A}^{\ast k}\otimes\mathcal{A}^{n}%
\longrightarrow a_{s^{\ast}}t\in\mathcal{A}^{n-k}.\label{act}%
\end{array}
\end{equation}
of annihilation operators for the algebraic description of our system. In this
way the triple $\{\Pi,P,a\}$, where $\Pi$ and $P$ are set of linear operators
and $a$ is the action of annihilation operators, is the initial data for our
model. The problem is to find and classify all triples of initial data which
lead to the well--defined models. The general solution for this problem is not
known for us. Hence we must restrict our attention for some examples.

If operators $P$ and $\Pi$ and the action $a$ of annihilation operators are
given in such a way that there is unique, nondegenerate, positive definite
scalar product, creation operators are adjoint to annihilation ones and vice
versa, then we say that we have a well--defined system with generalized statistics.

\emph{Example 1.} We assume here that $\Pi_{n}\equiv P_{n}\equiv
id_{E^{\otimes n}}$. This means that the algebra of states $\mathcal{A}$ is
identical with the full tensor algebra $TE$ over the space $E$, and the second
algebra $\mathcal{A}^{\ast}$ is identical with the tensor algebra $TE^{\ast}$.
The action (\ref{act}) of annihilation operators is given by the formula
$a_{x^{\ast i_{k}}\otimes\cdots\otimes x^{\ast i_{1}}}(x^{j_{1}}\otimes
\cdots\otimes x^{j_{n}}):=\delta_{i_{1}}^{j_{1}}\cdots\delta_{i_{k}}^{j_{k}%
}\ x^{j_{n-k+1}}\otimes\cdots\otimes x^{j_{n}}$. For the scalar product we
have the equation $\langle i_{n}\cdots i_{1}|j_{1}\cdots j_{n}\rangle
^{n}:=\delta^{i_{1}j_{1}}\cdots\delta^{i_{n}j_{n}}$. It is easy to see that we
have the relation and $a_{x^{\ast i}}a_{x_{j}}:=\delta_{i}^{j}\mathbf{1}$. In
this way we obtain the most simple example of well--defined system with
generalized statistics. The corresponding statistics is the so--called
infinite (Bolzman) statistics \cite{owg,gre}.

\emph{Example 2:} For this example we assume that $\Pi_{n}\equiv
id_{E^{\otimes n}}$. This means that $\mathcal{A}\equiv TE$ and $TE^{\ast}$.
For the scalar product and for the action of annihilation operators we assume
that there is a linear and invertible operator $T:E^{\ast}\otimes
E\longrightarrow E\otimes E^{\ast}$ defined by its matrix elements $T(x^{\ast
i}\otimes x^{j})=\Sigma_{k,\ast l}\ T_{k\ast l}^{\ast ij}\ x^{k}\otimes
x^{\ast l}$, such that we have $(T_{k\ast l}^{\ast ij})^{\ast}=\overline
{T}_{l\ast k}^{\ast ji}$, i.e. $T^{\ast}=\overline{T}^{t},$and $(T^{t})_{k\ast
l}^{\ast ij}=T_{l\ast k}^{\ast ji}$. Note that this operator not need to be
linear, one can also consider the case of nonlinear one. We also assume that
the operator $T^{\ast}$ act to the left, i.e. we have the relation
\begin{equation}
(x^{\ast j}\otimes x^{i})T^{\ast}=\Sigma_{l,\ast k}\ (x^{l}\otimes x^{\ast
k})\ \overline{T}_{l\ast k}^{\ast ji}, \label{teo1}%
\end{equation}
and $(T(x^{\ast i}\otimes x^{j}))^{\ast}\equiv(x^{\ast j}\otimes x^{i}%
)T^{\ast}$.The operator $T$ given by the formula (\ref{teo1}) is said to be
\textit{a twist} or \textit{a cross} operator. The operator $T$ describes the
cross statistics of quasiparticles and quasiholes. The set $P$ of projectors
is defined by induction $P_{n+1}:=(id\otimes P_{n})\circ R_{n+1}$, where
$P_{1}\equiv id$ and the operator $R_{n}$ is given by the formula
$R_{n}:=id+\tilde{T}^{(1)}+\tilde{T}^{(1)}\tilde{T}^{(2)}+\cdots+\tilde
{T}^{(1)}\dots\tilde{T}^{(n-1)}$, where $\tilde{T}^{(i)}:=id_{E}\otimes
\cdots\otimes\tilde{T}\otimes\cdots\otimes id_{E}$, $\tilde{T}$ on the $i$--th
place, and $(\tilde{T})_{kl}^{ij}=T_{l\ast j}^{\ast ki}$. If the operator
$\tilde{T}$ is a bounded operator acting on some Hilbert space such that we
have the following Yang-Baxter equation on $E\otimes E\otimes E$
\begin{equation}
(\tilde{T}\otimes id_{E})\circ(id_{E}\otimes\tilde{T})\circ(\tilde{T}\otimes
id_{E})=(id_{E}\otimes\tilde{T})\circ(\tilde{T}\otimes id_{E})\circ
(id_{E}\otimes\tilde{T}),
\end{equation}
and $||\tilde{T},||\leq1$, then according to Bo$\dot{z}$ejko and Speicher
\cite{bs2} there is a positive definite scalar product
\begin{equation}
\langle s|t\rangle_{T}^{n}:=\langle s|P_{n}t\rangle_{0}^{n} \label{csca}%
\end{equation}
for $s,t\in\mathcal{A}^{n}\equiv E^{\otimes n}$. Note that the existence of
nontrivial kernel of operator $P_{2}\equiv R_{1}\equiv id_{E\otimes E}%
+\tilde{T}$ is essential for the nondegeneracy of the scalar product
\cite{jswe}. One can see that if this kernel is trivial, then we obtain
well--defined system with generalized statistics \cite{m10,ral}.

\emph{Example 3:} If the kernel of $P_{2}$ is nontrivial, then the scalar
product (\ref{csca}) is degenerate. Hence we must remove this degeneracy by
factoring the mentioned above scalar product by the kernel. We assume that
there is an ideal $I\subset TE$ generated by a subspace $I_{2}\subset
\mathrm{ker}P_{2}\subset E\otimes E$ such that $a_{s^{\ast}}I\subset I$ for
every $s^{\ast}\in\mathcal{A}^{\ast}$, and for the corresponding ideal
$I^{\ast}\subset E^{\ast}\otimes E^{\ast}$ we have $a_{s^{\ast}}t=0$ for every
$t\in TE$ and $s^{\ast}\in I^{\ast}$. The above ideal $I$ is said to be Wick
ideal \cite{jswe}. We have here the following formulae $\mathcal{A}%
:=TE/I,\quad\mathcal{A}^{\ast}:=TE^{\ast}$ for our algebras. The projection
$\Pi$ is the quotient map $\Pi:\tilde{s},\in TE\longrightarrow s\in
TE/I\equiv\mathcal{A}$. For the scalar product we have here the following
relation $\langle s|t\rangle_{B,T}:=\langle\tilde{s}|\tilde{t}\rangle_{T}$ for
$s=P_{m}(\tilde{s})$ and $t=P_{n}(\tilde{t})$. One can define here the action
of annihilation operators in such a way that we obtain well--defined system
with generalized statistics \cite{ral}.

\emph{Example 4:} If $B=\frac{1}{\mu}\tilde{T},$, where $\mu$ is a parameter,
then the third condition (\ref{cd}) is equivalent to the well known Hecke
condition for $\tilde{T},$ and we obtain the well--known relations for Hecke
symmetry and quantum groups \cite{P,PW,Ke}.

\textbf{Physical effect}. Let us consider the system equipped with generalized
statistics and described by two operators $T$ and $B$ like in Example 4. We
assume here in addition that a linear and Hermitian operator $S:E\otimes
E\longrightarrow E\otimes E$ such that
\begin{equation}
S^{\left(  1\right)  }S^{\left(  2\right)  }S^{\left(  1\right)  }=S^{\left(
2\right)  }S^{\left(  1\right)  }S^{\left(  2\right)  },\quad\mbox{and}\quad
S^{2}=id_{E\otimes E}%
\end{equation}
is given. If we have the following relation $\tilde{T},\equiv B\equiv S$, then
it is easy to see that the conditions (\ref{cd}) are satisfied and we have
well--defined system with generalized statistics. Let us assume for simplicity
that the operator $S$ is diagonal and is given by the following equation
$S(x^{i}\otimes x^{j})=\epsilon^{ij}x^{j}\otimes x^{i}$, for $i,j=1,\ldots,N$,
where $\epsilon^{ij}\in\mathbb{C}$, and $\epsilon^{ij}\epsilon^{ji}=1$. In
general we have $c_{ij}=-(-1)^{\Sigma_{ij}}q^{\Omega_{ij}}$, where
$\Sigma:=(\Sigma_{ij})$ and $\Omega:=(\Omega_{ij})$ are integer--valued
matrices such that $\Sigma_{ij}=\Sigma_{ji}$ and $\Omega_{ij}=-\Omega_{ji}$.
The algebra $\mathcal{A}$ is here a quadratic algebra generated by relations
\begin{equation}
x^{i}\odot x^{j}=\epsilon^{ij}x^{j}\odot x^{i},\quad\mbox{and}\quad(x^{i}%
)^{2}=0\quad\mbox{if}\quad\epsilon^{ii}=-1
\end{equation}
We also assume that $\epsilon^{ii}=-1$ for every $i=1,\ldots,N$. In this case
the algebra $\mathcal{A}$ is denoted by $\Lambda_{\epsilon}(N)$. Now let us
study the algebra $\Lambda_{\epsilon}(2)$, where $\epsilon^{ii}=-1$ for
$i=1,2$, and $\epsilon^{ij}=1$ for $i\neq j$, in more details. In this case
our algebra is generated by $x^{1}$ and $x^{2}$ such that we have $x^{1}\odot
x^{2}=x^{2}\odot x^{1},\quad(x^{1})^{2}=(x^{2})^{2}=0$. Note that the algebra
$\Lambda_{\epsilon}(2)$ is an example of the so--called $Z_{2}\oplus Z_{2}%
$--graded commutative color Lie superalgebra \cite{luri}. Such algebra can be
transformed into the usual Grassmann algebra $\Lambda_{2}$ generated by
$\Theta^{1}$ and $\Theta^{2}$ such that we have the anticommutation relation
$\Theta^{1}\ \Theta^{2}=-\Theta^{2}\ \Theta^{1}$, and $(\Theta^{1}%
)^{2}=(\Theta^{1})^{2}=0$. In order to do such transformation we use the
Clifford algebra $C_{2}$ generated by $e^{1},e^{2}$ such that we have the
relations $e^{i}\ e^{j}+e^{j}\ e^{i}=2\delta^{ij}$, for $i,j=1,2$. For
generators $x^{1}$, and $x^{2}$ of the algebra $\Lambda_{\epsilon}(2)$ the
transformation is given by $\Theta^{1}:=x^{1}\otimes e^{1}$, for $\Theta
^{2}:=x^{2}\otimes e^{2}$. It is interesting that the algebra $\Lambda
_{\epsilon}(2)$ can be represented by one Grassmann variable $\Theta$,
$\Theta^{2}=0$ as follows $x^{1}=(\Theta,1),\quad x^{2}=(1,\Theta)$. For the
product $x^{1}\odot x^{2}$ we obtain $x^{1}\odot x^{2}=(\Theta,\Theta)$.

In physical interpretation generators $\Theta^{1}$ and $\Theta^{2}$ of the
algebra $\Lambda_{2}$ represents two fermions. They anticommute and according
to the Pauli exclusion principle we can not put them into one energy level.
Observe that the corresponding generators $x^{1}$ and $x^{2}$ of the algebra
$\Lambda_{\epsilon}(2)$ commute, their squares disappear and they describe two
different quasiparticles. This means that these quasiparticles behave
partially like bosons, we can put them simultaneously into one energy levels,
and single fermion can be transformed under certain interactions into a system
of two different quasiparticles.

\section{Interaction Processes and Entwining Structures
\footnotep{Last change October 19, 2002}}

Let us assume that an interacting system of charged particles and quanta is
given as a starting point for the study of particle processes. Our main
assumption is that there is a set of initial configurations of the system and
there is also a set of final configurations representing possible results of
interactions. An arbitrary initial configuration can be transform into final
one as a result of interaction processes. Every such transformation is said to
be a \textit{process}. Such transformation can be done in several different
way. All of them should produce the same result. Hence there is the uniqueness
problem for such description. If all final configurations for the system under
considerations can be described in an unique way as a result of transformation
of an initial configuration, then we say that the system is equipped with a
\textit{coherent process}.

In this paper possible particle processes are studied in terms of monoidal
categories. Our considerations are based on assumptions discussed in previous
publications \cite{top,castat,sin,wmq} and shortly summarized in the
Preliminaries. Particle systems equipped with possibility to absorption or
emission of certain quanta are described in terms of a category of modules and
comodules. The process of absorption of quanta of the external field can be
described as a left or right coaction of $\mathcal{C}$ on $\mathcal{M}$.
Dually, the process of emission of quanta can be described as a left or right
action of $\mathcal{C}$ on $\mathcal{M}$. In this case the coherence means
that left and right (co-)action leads to the equivalent result. In our
approach the equivalence can be expressed by an invertible entwining map
$\tau:\mathcal{C}\otimes\mathcal{A}\rightarrow\mathcal{A}\otimes\mathcal{C}$
\cite{brmj}. Composite systems which contain both charged particles and quanta
are described as entwined modules. It is shown that under some assumptions
there is an algebra $\mathcal{W}(\mathcal{A},\mathcal{C},\tau)$ whose
representations represent the corresponding quantum states.

Let $\mathfrak{C}=\mathfrak{C}(\otimes,k)$ be a monoidal category with duals.
The monoidal operation $\otimes:\mathfrak{C}\times\mathfrak{C}\rightarrow
\mathfrak{C}$ is a bifunctor which has a two-sided identity object $k$. The
category $\mathfrak{C}$ can contain some special objects like algebras,
coalgebras, modules or comodules, etc... According to our previous
considerations all possible particle processes can be represented as arrows of
the category $\mathfrak{C}$. In our case $k=\mathbb{C}$ is the field of
complex numbers. If $f:\mathcal{U}\rightarrow\mathcal{V}$ is an arrow from
$\mathcal{U}$ to $\mathcal{V}$, then the object $\mathcal{U}$ represents
physical objects before interactions and $\mathcal{V}$ represents possible
results of interactions. We assume that different objects of the underlying
category describe physical objects of different nature, charged particles,
quasiparticles or different species of quanta of an external field, etc... Let
$\mathcal{U}$ be an object of the category $\mathfrak{C}$, then the object
$\mathcal{U}^{\ast}$ corresponds for antiparticles, holes or quasiholes or
dual field, respectively. If $\mathcal{U}$ and $\mathcal{V}$ are two different
objects of the category $\mathcal{C}$, then the product $\mathcal{U}%
\otimes\mathcal{V}$ is also an object of the category, it represents a
composite system composed from object of different nature.

Let $\mathcal{A}$ be an unital associative algebra, and $\mathcal{C}$ be a
counital coassociative coalgebra in $\mathfrak{C}$. We use the notation
$\Delta(f):=f^{(1)}\otimes f^{(2)}$ for the coalgebra comultiplication. In our
approach charged particles are represented by the algebra $\mathcal{A}$, an
external quantum field is characterized by coalgebra $\mathcal{C}$. Composite
systems which contain both charged particles and quanta can be described as a
product of copies of $\mathcal{A}$ and $\mathcal{C}$. The multiplication
$m:\mathcal{A}\otimes\mathcal{A}\rightarrow\mathcal{A}$ represents the
creation process of a single classical object from a composite system of
objects of the same species. The result of comultiplication $\triangle
:\mathcal{B}\rightarrow\mathcal{B}\otimes\mathcal{B}$ represents a
'composition' process from copies of objects of the same nature. Let
$\mathcal{C}$ be a coalgebra and $\mathcal{C}^{\ast}$ be an algebra in
duality, i. e. we have a bilinear pairing $\langle-,-\rangle:\mathcal{C}%
\otimes\mathcal{C}^{\ast}\rightarrow k$ such that $<\triangle f,s\otimes
t>=<f,m(s\otimes t)>$, where $f\in\mathcal{C}$ and $s,t\in\mathcal{C}^{\ast}$.
In our physical interpretation this duality indicates certain relations
between processes.

\textbf{Entwined structures}. Let $\mathcal{A}$ be an algebra, and
$\mathcal{C}$ be a coalgebra. We denote by $M_{\mathcal{A}}^{\mathcal{C}}$ the
category of right $\mathcal{A}$-modules and right $\mathcal{C}$-comodules. A
mapping $\tau:\mathcal{C}\otimes\mathcal{A}\rightarrow\mathcal{A}%
\otimes\mathcal{C}$ such that
\begin{equation}%
\begin{array}
[c]{ll}%
\tau\circ(id\otimes m)=(m\otimes id)\circ\tau_{23}\circ\tau_{12}, & \tau
\circ(f\otimes1)=1\otimes f,\\
(id\otimes\Delta)\circ\tau=\tau_{12}\circ\tau_{23}\circ(\Delta\otimes id), &
(id\otimes\varepsilon)\circ\tau=\varepsilon\otimes id,
\end{array}
\end{equation}
is said to be entwining, where $\tau_{12}:=\tau\otimes id$, and $f\in
\mathcal{C}$. Let $\mathcal{A}$ be a $\mathcal{C}$-Galois extension of
$\mathcal{B}$, then there is a unique entwining map $\tau:\mathcal{C}%
\otimes\mathcal{A}\rightarrow\mathcal{A}\otimes\mathcal{C}$ (see \cite{brmj}).
The entwining map is given by $\tau(f\otimes a):=\beta\circ(id\otimes
_{\mathcal{B}}m)\circ(\chi\otimes id)$, where $\chi:=\beta^{-1}(1\otimes
f),f\in\mathcal{C}$. We use the following notation $\tau(f\otimes
a):=a_{(1)}\otimes f_{(2)}$ for the entwining $\tau$ and $\mathcal{A}%
^{n}:=\underbrace{\mathcal{A}\otimes\cdots\otimes\mathcal{A}}_{n}$for tensor
product of $n$ copies of $\mathcal{A}$.

Let us define a mapping $\Psi_{1n}:\mathcal{C}\otimes\mathcal{A}%
^{n}\rightarrow\mathcal{A}^{n}\otimes\mathcal{C}$ by the relations $\Psi
_{11}:=\tau,\quad\Psi_{1n}:=\underbrace{(id\otimes\tau)\circ\cdots\circ
(\tau\otimes id)}_{n}$. If $\tau$ is invertible, then we also define
$\Psi_{n1}:=(\Psi_{1n})^{-1}$.

\textbf{Entwined modules}. A right $\mathcal{A}$-module $\mathcal{M}$ equipped
with an action $\alpha:\mathcal{M}\otimes\mathcal{A}\rightarrow\mathcal{M}$
which is also a right $\mathcal{C-}$ comodule with a coaction $\delta
:\mathcal{M}\rightarrow\mathcal{M}\otimes\mathcal{C}$ such that
\begin{equation}%
\begin{array}
[c]{ccccccc}%
\mathcal{M}\otimes\mathcal{A} & \rightarrow & \rightarrow & \mathcal{M} &
\rightarrow & \rightarrow & \mathcal{M}\otimes\mathcal{C}\\
\parallel &  &  &  &  &  & \parallel\\
\mathcal{M}\otimes\mathcal{A} & \overset{\delta}{\rightarrow} & \mathcal{M}%
\otimes\mathcal{C}\otimes\mathcal{A} & \overset{id\otimes\tau}{\rightarrow} &
\mathcal{M}\otimes\mathcal{A}\otimes\mathcal{C} & \overset{\alpha}%
{\rightarrow} & \mathcal{M}\otimes\mathcal{C}%
\end{array}
\end{equation}
is said to be an entwined module or $(\mathcal{A},\mathcal{C},\tau)$-module.
Let $\mathcal{M}$ and $\mathcal{N}$ be two $(\mathcal{A},\mathcal{C},\tau
)$-modules. A mapping $f:\mathcal{M}\rightarrow\mathcal{N}$ such that
\begin{equation}%
\begin{array}
[c]{rcclrcl}%
\mathcal{A}\otimes\mathcal{M} & \overset{\alpha}{\rightarrow} & \mathcal{M} &
& \mathcal{M} & \overset{\delta}{\rightarrow} & \mathcal{M}\otimes
\mathcal{C}\\
id\otimes f\downarrow &  & f\downarrow &  & \downarrow f &  & \downarrow
f\otimes id\\
\mathcal{A}\otimes\mathcal{N} & \overset{\alpha}{\rightarrow} & \mathcal{N} &
& \mathcal{N} & \overset{\delta}{\rightarrow} & \mathcal{N}\otimes\mathcal{C}%
\end{array}
\end{equation}
is said to be a $(\mathcal{A},\mathcal{C},c)$-module morphism. It is obvious
that the collection of all $(\mathcal{A},\mathcal{C},c)$-modules and
$(\mathcal{A},\mathcal{C},c)$-module morphisms is a category. We denote this
category by $\mathcal{M}(\mathcal{A},\mathcal{C},c)$. For an algebra
$\mathcal{A}$, and a coalgebra $\mathcal{C}$ with a coaction $\delta
:\mathcal{A}\rightarrow\mathcal{A}\otimes\mathcal{C}$ we define $\mathcal{C}%
_{\mathcal{A}}^{n}:=\underbrace{\mathcal{C}\otimes\cdots\otimes\mathcal{C}%
}_{n}\otimes\mathcal{A}$, $n=1,2,\ldots.$ The space $\mathcal{C}_{\mathcal{A}%
}^{n}$ is (i) a right $\mathcal{A}$-module, (ii) a right $\mathcal{C}%
$-comodule, and (iii) a $(\mathcal{A},\mathcal{C},\tau)$-module,
$n=1,2,\ldots$.

We define the action $\alpha_{n}:\mathcal{C}_{\mathcal{A}}^{n}\otimes
\mathcal{A}\rightarrow\mathcal{A}$ and the coaction $\delta_{n}:\mathcal{C}%
_{\mathcal{A}}^{n}\rightarrow\mathcal{C}_{\mathcal{A}}^{n}\otimes\mathcal{C}$
by the following formulae $\alpha_{n}:=id\otimes m,\quad\delta_{n}%
:=id\otimes\delta$. This simply follows from the associativity and coassociativity.

\textbf{Crossed product}. Let $\mathfrak{C}=(\mathfrak{C}_{0},\mathfrak{C}%
_{1})$ be a category whose objects $\mathfrak{C}_{0}$ are associative and
unital algebras over a field $k$ and whose morphisms are algebra morphisms.
For our purposes here, we shall denote by $\mathcal{A}\otimes\mathcal{B}$ the
tensor product $\mathcal{A}\otimes_{k}\mathcal{B}$, where $\mathcal{A}%
,\mathcal{B}$ are considered as $k$-linear spaces. A linear mapping
$\Psi:\mathcal{B}\otimes\mathcal{A}\rightarrow\mathcal{A}\otimes\mathcal{B}$
such that we have the following relations
\begin{align}
\Psi\circ(id_{\mathcal{B}}\otimes m_{\mathcal{A}})&=(m_{\mathcal{A}}\otimes
id_{\mathcal{B}})\circ(id_{\mathcal{A}}\otimes\Psi)\circ(\Psi\otimes
id_{\mathcal{A}}),\\ \Psi\circ(m_{\mathcal{B}}\otimes id_{\mathcal{A}%
})&=(id_{\mathcal{A}}\otimes m_{\mathcal{B}})\circ(\Psi\otimes id_{\mathcal{B}%
})\circ(id_{\mathcal{B}}\otimes\Psi)
\end{align}
is said to be an algebra cross or twist \cite{bma}. We use here the notation
$\Psi(b\otimes a)=\Sigma a_{(1)}\otimes b_{(2)}$ for $a\in\mathcal{A}%
,b\in\mathcal{B}$. The tensor product $\mathcal{A}\otimes_{k}\mathcal{B}$ of
algebras $\mathcal{A}$ and $\mathcal{B}$ equipped with the multiplication
\begin{equation}%
\begin{array}
[c]{c}%
m_{\Psi}:=(m_{\mathcal{A}}\otimes m_{\mathcal{B}})\circ(id_{\mathcal{A}%
}\otimes\Psi\otimes id_{\mathcal{B}})
\end{array}
\end{equation}
is an associative algebra called a crossed product with respect to the cross
symmetry $\Psi$ \cite{bma} and it is denoted by $\mathcal{W}=\mathcal{W}%
_{\Psi}(\mathcal{A},\mathcal{B})=\mathcal{A}\otimes_{\Psi}\mathcal{B}$. There
is one to one correspondence between algebra cross $\Psi:\mathcal{B}%
\otimes\mathcal{A}\rightarrow\mathcal{A}\otimes\mathcal{B}$ and crossed
product $\mathcal{W}$ of algebras $\mathcal{A}$ and $\mathcal{B}$. Let
$(\mathcal{A},\mathcal{C},\tau)$ an entwining structure $(\mathcal{A}%
,\mathcal{C},\tau)$, and $\mathcal{C}^{\ast}:=\hom(\mathcal{C},k)$, be the
dual of $\mathcal{C}$. Observe that there is the bilinear pairing
$\langle-,-\rangle:\mathcal{C}\otimes\mathcal{C}^{\ast}\rightarrow k$ defined
by $\langle f,x\rangle\equiv ev(f\otimes x):=f(x)$. If there exists a unique
mapping $\widetilde{\tau}:\mathcal{A}\otimes\mathcal{C}^{\ast}\rightarrow
\mathcal{C}^{\ast}\otimes\mathcal{A}$ such that
\begin{equation}%
\begin{array}
[c]{ccl}%
\mathcal{C}\otimes\mathcal{A}\otimes\mathcal{C}^{\ast} & \overset
{id\otimes\widetilde{\tau}}{\longrightarrow} & \mathcal{C}\otimes
\mathcal{C}^{\ast}\otimes\mathcal{A}\\
\downarrow\tau\otimes id &  & \downarrow ev\otimes id\\
\mathcal{A}\otimes\mathcal{C}\otimes\mathcal{C}^{\ast} & \overset{id\otimes
ev}{\longrightarrow} & \mathcal{A},
\end{array}
\end{equation}
then the entwining structure $(\mathcal{A},\mathcal{C},\tau)$ is said to be
$\mathcal{C}^{\ast}$-factorizable. If $(\mathcal{A},\mathcal{C},\tau)$ a
factorizable entwining structure, then the mapping $\widetilde{\tau
}:\mathcal{A}\otimes\mathcal{C}^{\ast}\rightarrow\mathcal{C}^{\ast}%
\otimes\mathcal{A}$ is an algebra cross.

Let $(\mathcal{A},\mathcal{C},\tau)$ a factorizable entwining structure
$(\mathcal{A},\mathcal{C},\tau)$, then there is a cross symmetry $\Psi$ and
the corresponding algebra cross tensor product $\mathcal{W}(\mathcal{A}%
,\mathcal{C},\tau)=\mathcal{A}\otimes_{\Psi}\mathcal{C}^{\ast}$ of algebras
$\mathcal{A}$ and $\mathcal{C}^{\ast}$. If $(\mathcal{A},\mathcal{C},\tau)$ a
factorizable entwining structure, then according to the last lemma there
exists a unique algebra cross $\widetilde{\tau}:\mathcal{A}\otimes
\mathcal{C}^{\ast}\rightarrow\mathcal{C}^{\ast}\otimes\mathcal{A}$.

\section{Particles interactions and quantum logic
\footnotep{Last change January 28, 2003}}

It is well-known that from algebraic point of view that $q$-deformed
commutation relations for particles equipped with generalized statistics can
be described in terms of the so-called Wick algebras \cite{jswe}. The
construction Wick algebras is based on a single operator $T$ which describe
the deformation. A second operator $B$ is required for consistency relations.
A proposal for the general algebraic formalism for description of particle
systems equipped with an arbitrary generalized statistics based on the concept
of monoidal categories with duality has been given by the Author in
\cite{castat}. The physical interpretation for this formalism was shortly
indicated. A few examples of applications for this formalism are considered in
\cite{top,com}.

In this paper we are going to study possible states of an interacting particle
systems in terms of quantum logics. Our notion of quantum logics is specific
and is partially based on the book \cite{ppp}. All our considerations are
based on previously developed algebraic formalism for particles with
generalized statistics \cite{gqsi,qcc}. We describe a system with generalized
statistics as a quantum logics representation. Our study is motivated by
possible physical applications to describe pseudoparticle configurations of
magnetic crystal and related topics, see \cite{lul} and reference therein.

\textbf{Discrete system and quantum logics}. The starting point for our study
is an interacting particle system. It is natural to assume that the whole
physical world is divided into two parts: a classical system and its quantum
environment. The classical system represents a physical reality, all that can
be observed. The quantum environment represents all quantum possibilities to
become a part of the reality in the future. Our fundamental assumption is that
an initial particle configuration is transformed under interaction into a
composite system consisting all possible results of interactions. We assume
that all possible composite systems as a results of interactions can be
constructed from an initial set of elementary ones. We would like to construct
an algebraic model representing these possibilities \cite{qcc}.

A system which contains a charge and certain number of quanta as a result of
interaction with the quantum field is said to be a \textit{dressed particle}
or a \textit{lattice}. We describe a dressed particle as a non-local discrete
system which contains $n$ centers (vertexes or lattice sites). Every center
can be dressed by quanta of $N$ different sorts. A center with a place for
certain quanta is said to be a \textit{quasihole} or \textit{an empty quantum
level}. A centre dressed with a single quantum of the field is said to be a
\textit{quasiparticle}.

Let us denote by $L_{n}:=(I_{n},\mathbb{Q})$ a collection of all possible
states corresponding to a dressed particle with $n$-centers and $N$-sorts
quanta. We assume that this collection contains a finite collection of
elementary (or primary) states $\mathbb{Q}\equiv\mathbb{Q}_{N}:=\{1,2,\ldots
,N<\infty\}$, a corresponding collection of $\ast$-conjugated states
$\mathbb{Q}^{\ast}\equiv\mathbb{Q}_{N}^{\ast}:=\{\ast N,\ldots,\ast2,\ast1\}$
and a collection of maps $L_{n}:=\{\sigma:k\in I_{n}\mapsto\sigma(k):=i_{k}%
\in\mathbb{Q}_{N}\cup\mathbb{Q}_{N}^{\ast}\cup\emptyset\}$, where
$k=1,\ldots,n$, $I_{n}:=\{1,\ldots,n\}$; $n=1,2,\dots$, called composite
states (or chain of states). This means that $\sigma\in L_{n}$ is a finite
sequence of length $n$ .i.e $\sigma:=\{i_{1},\ldots,i_{n}\}$, $i_{k}%
\in\mathbb{Q}_{N}\cup\mathbb{Q}_{N}^{\ast}\cup\emptyset,k=1,\ldots,n$. The
empty sequence is denoted by $\sigma(0):=\emptyset$. We assume that
$L_{0}:=\emptyset$ for $n=0$. There is a series of actions $\varrho_{n}:\pi\in
S_{n}\rightarrow\varrho_{n}(\pi)\in End(L_{n})$, where $S_{n}$ is the $n$-th
symmetric group, and
\begin{equation}
\varrho(\pi):=\left(
\begin{array}
[c]{l}%
\sigma\\
\sigma\circ\pi^{-1}%
\end{array}
\right)  ,\;\mbox{}\;\sigma\in L,\;\pi\in S_{n}.
\end{equation}

For $\sigma:=\{i_{1},\ldots,i_{n}\}$ we define an $\ast$-operation by the
formula $\ast\sigma:=\{\ast i_{n},\ldots,\ast i_{1}\}$, where $\ast(i):=\ast
i$, $\ast(\ast i):=\ast\ast i:=i$. Observe that 1) $\emptyset\in L$; 2) if
$\sigma\in L_{n}$, then there is $\ast\sigma\in L_{n}$; 3) if $\sigma_{n}\in
L_{n}$ for $n=1,2,\ldots$, then $\bigcup_{n=1,2,\ldots}\sigma_{n}\in L$; and
$\ast(\ast\sigma)=\sigma,\;\sigma\cap\ast\sigma=\emptyset$. We say that
$L:=\bigcup_{n=0,1,\ldots}\;L_{n}$ is a family of quantum logics if and only
if an transitive and anti-reflexive relation $\perp$ in $L$ is defined. The
relation $\perp$ is said to be an orthogonality, we denote by $\sigma
\perp\sigma^{\prime}$ a pair of orthogonal elements of $L$ \cite{ppp}. We
assume in addition that $\{\ast i\}\perp\{j\}$, if $i\neq j$.

Let $\sigma,\sigma^{\prime}\in L$ and $\sigma\perp\sigma^{\prime}$, then the
pair $(\sigma,\sigma^{\prime})$ is said to be a mutually excluded. A
collection $\{\sigma_{i_{1}},\ldots,\sigma_{i_{n}}\}$ is said to a complete if
and only if $\sigma_{i_{1}}\perp\ldots\perp\sigma_{i_{n}}=\mathbb{Q}%
\cup\mathbb{Q}^{\ast}$. Obviously $\mathbb{Q}\cup\mathbb{Q}^{\ast}$ is
complete $\{i_{1}\}\perp\ldots\perp\{i_{N}\}=\mathbb{Q}\cup\mathbb{Q}%
_{N}^{\ast}$.

Let $L$ and $L^{\prime}$ be quantum logics. A mapping $f:L\rightarrow
L^{\prime}$ such that 1) $f(\emptyset)=\emptyset$; 2) $f(\ast\sigma)=\ast
f(\sigma)$; 3) $\mbox{if}\;\sigma_{i_{1}}\perp\cdots\perp\sigma_{i_{n}%
},\mbox{then}\;f(\sigma_{i_{1}})\perp\cdots\perp f(\sigma_{i_{n}})$ is said to
be a logic morphism.

\textbf{Linear representations}. A (symmetric) operad $E$ is defined as a
collection of sets $\mathbf{E}:=\{E(k):k=0,1,2,\ldots,n,\ldots\}$, equipped
with a collection of structure mappings
\begin{equation}%
\begin{array}
[c]{c}%
\gamma_{k_{1}\cdots k_{n}}:E(n)\times E(k_{1})\times\cdots\times
E(k_{n})\rightarrow E(k_{1}+\cdots+k_{n})
\end{array}
\end{equation}
for every $k_{1},\ldots,k_{n}=1,2,\ldots$, $n=1,2,\ldots$ satisfying some
known compatibility condition for composition and symmetric group actions
\cite{opm}. We use the notation $\gamma_{k_{1}\cdots k_{n}}(v;v_{1}%
,\cdots,v_{n}):=v(v_{1},\cdots,v_{n})$, where $v\in E(n),v_{1}\in
E(k_{1}),\cdots,v_{n}\in E(k_{n})$. We assume that there is an element
$\mathbf{1}\in E(0)$ called the unit.

Let $L:=(\mathbb{Q},P)$ be a quantum logic. A representation of $L$ In an
operad $\mathbf{E}$ is a map $x:L\rightarrow\mathbf{E}$ such that
\begin{equation}%
\begin{array}
[c]{ll}%
x(\sigma)\in E(n), & \mbox{for}\;\sigma:=\{i_{1},\ldots,i_{n}\}\\
x(\emptyset):=\mathbf{1}, & x(\ast\sigma)=(x(\sigma))^{\ast};
\end{array}
\end{equation}
and $x(\sigma_{i_{1}})\perp\cdots\perp x(\sigma_{i_{n}})$ for $\sigma_{i_{1}%
}\perp\cdots\perp\sigma_{i_{n}}$. For the action $\varrho_{n}$ of the
symmetric group we assume that $x(\varrho_{n}(\pi)\sigma)=\varrho_{n}^{E}\circ
x(\sigma)$, where $\varrho^{E}$ is the corresponding action in the operad.

Here we restrict our attention to the particular choice $E(n):=E^{\otimes n}$,
$E(0):=\mathbf{1}I$, $E(1):=E$, where $E$ is a linear space over a field $I$
equipped with a basis $S:=\{x(\{i\}):=x^{i}:i=1,\ldots,N\}$. We have
\begin{equation}
x(\bigcup_{i=1,2,\ldots}\sigma_{i}):=\bigotimes_{i=1,2,\ldots}x(\sigma_{i}).
\end{equation}
For $\sigma:=\{i_{1},\ldots,i_{n}\}\in P$ we have an extension $x(\sigma
):=x^{i_{1}}\otimes\cdots\otimes x^{i_{n}}$. Similarly $S^{\ast}%
:=\{x(\ast\sigma_{i}):=x^{\ast i}:i=1,\ldots,N\}$ forms a basis of the
conjugate space $E^{\ast}$. The pairing $g_{E}:E^{\ast}\otimes E\rightarrow I$
and the corresponding scalar product is given by $g_{E}(x^{\ast i}\otimes
x^{j})\equiv(x^{\ast i}|x^{j})=\langle x^{i}|x^{j}\rangle:=\delta^{ij}$.

Let $T:E^{\ast}\otimes E\rightarrow E\otimes E^{\ast}$ be a linear and
Hermitian operator with matrix elements
\begin{equation}
T(x^{\ast i}\otimes x^{j})=\Sigma\ T_{kl}^{ij}\;x^{k}\otimes x^{\ast l},
\label{cross}%
\end{equation}
then the quotient $\mathcal{W}(T)=T(E\oplus E^{\ast})/I_{T}$, where the ideal
$I_{T}$ is given by the relation
\begin{equation}
I_{T}:=gen\{x^{\ast i}\otimes x^{j}-\Sigma\;T_{kl}^{ij}\;x^{k}\otimes x^{\star
l}-(x^{\star i}|x^{j})\}
\end{equation}
is said to be Hermitian Wick algebra \cite{jswe}.

\textbf{Algebra representation and cross product}. Let $L$ be a quantum logic,
$\mathcal{A}$ be an an unital and associative algebra,
\begin{equation}%
\begin{array}
[c]{l}%
\mathcal{A}:=\bigoplus\limits_{n}\;\mathcal{A}^{n},
\end{array}
\end{equation}
equipped with an associative multiplication $m:\mathcal{A}\otimes
\mathcal{A}\rightarrow\mathcal{A}$ then there is a representation
$x:L\rightarrow\mathcal{A}$ of $L$,
such that $x(\sigma):=m_{n}(x^{i_{1}}\otimes\cdots\otimes x^{i_{n}})$ for
$\sigma:=\{i_{1},\ldots,i_{n}\}\in P$, and $m_{n}:=\underbrace{m(id\otimes
\cdots(id\otimes m)\cdots)}_{n}$. A pair of algebras $\mathcal{A}$ and
$\mathcal{A}^{\ast}$ is said to be conjugated algebras if and only if there is
an antilinear and involutive isomorphism $(-)^{\ast}:\mathcal{A}%
\rightarrow\mathcal{A}^{\ast}$, i. e. we have the relations $m_{\mathcal{A}%
^{\ast}}(b^{\ast}\otimes a^{\ast})=(m_{\mathcal{A}}(a\otimes b))^{\ast}%
,\quad(a^{\ast})^{\ast}=a$, where $a,b\in\mathcal{A}$ and $a^{\ast},b^{\ast}$
are their images under the isomorphism $(-)^{\ast}$.

Let $(\mathcal{A},\mathcal{A}^{\ast})$ be a pair of conjugate algebras. A
linear mapping $\Psi:\mathcal{A}^{\ast}\otimes\mathcal{A}\rightarrow
\mathcal{A}\otimes\mathcal{A}^{\ast}$ such that $\Psi|_{E^{\ast}\otimes
E}=T+g_{E}$, and we have the following relations \cite{csv,bma}
\begin{equation}%
\begin{array}
[c]{l}%
\Psi\circ(id_{\mathcal{A}^{\ast}}\otimes m_{\mathcal{A}})=(m_{\mathcal{A}%
}\otimes id_{\mathcal{A}^{\ast}})\circ(id_{\mathcal{A}}\otimes\Psi)\circ
(\Psi\otimes id_{\mathcal{A}}),\\
\Psi\circ(m_{\mathcal{A}^{\ast}}\otimes id_{\mathcal{A}})=(id_{\mathcal{A}%
}\otimes m_{\mathcal{A}^{\ast}})\circ(\Psi\otimes id_{\mathcal{A}^{\ast}%
})\circ(id_{\mathcal{A}^{\ast}}\otimes\Psi)\label{twc}%
\end{array}
\end{equation}
is said to be a cross symmetry or generalized braidings generated by $T$. We
use here the notation $\Psi(b^{\ast}\otimes a)=\Sigma a_{(1)}\otimes
b_{(2)}^{\ast}$ for $a\in\mathcal{A},b^{\ast}\in\mathcal{A}^{\ast}$. The
tensor product $\mathcal{A}\otimes\mathcal{A}^{\ast}$ equipped with the
multiplication
\begin{equation}
m_{\Psi}:=(m_{\mathcal{A}}\otimes m_{\mathcal{A}^{\ast}})\circ(id_{\mathcal{A}%
}\otimes\Psi\otimes id_{\mathcal{A}^{\ast}}) \label{mul}%
\end{equation}
is an associative algebra isomorphic to the Hermitian Wick algebra \cite{bma}
and it is denoted by $\mathcal{W}=\mathcal{W}_{\Psi}(\mathcal{A}%
)=\mathcal{A}>\!\!\!\lhd_{\Psi},\mathcal{A}^{\ast}$. Let $H$ be a linear
space. We denote by $L(H)$ the algebra of linear operators acting on $H$.

Let $\mathcal{W}\equiv\mathcal{A}>\!\!\!\lhd_{\Psi},\mathcal{A}^{\ast}$ be a
Hermitian Wick algebra. If $\pi_{\mathcal{A}}:\mathcal{A}\rightarrow L(H)$ is
a representation of the algebra $\mathcal{A}$, such that we have the relation
\begin{equation}%
\begin{array}
[c]{l}%
(\pi_{\mathcal{A}}(b))^{\ast}\pi_{\mathcal{A}}(a)=\Sigma\pi_{\mathcal{A}%
}(a_{(1)})\pi_{\mathcal{A}^{\ast}}(b_{(2)}^{\ast}),\\
\pi_{\mathcal{A}^{\ast}}(a^{\ast}):=(\pi_{\mathcal{A}}(a))^{\ast},\label{wre}%
\end{array}
\end{equation}
then there is a representation $\pi_{\mathcal{W}}:\mathcal{W}\rightarrow L(H)$
of the algebra $\mathcal{W}$ \cite{bma}. The relations (\ref{wre}) are said to
be a commutation relation if there is a positive definite scalar product on
$\mathcal{A}$ such that $\langle\pi_{\mathcal{A}^{\ast}}(x^{\ast}%
)f|g\rangle=\langle f|\pi_{\mathcal{A}}(x)g\rangle$. Note that if we use the
notation $\pi_{\mathcal{A}}(x^{i})\equiv a_{x^{i}}^{+},\quad\pi_{\mathcal{A}%
^{\ast}}(x^{\ast i})\equiv a_{x^{\ast i}}$, and the cross $T$ is given by its
matrix elements (\ref{cross}), then the commutation relations (\ref{wre}) can
be given in the following form
\begin{equation}
a_{x^{\ast i}}a_{x^{j}}^{+}-T_{kl}^{ij}\;a_{x^{l}}^{+}a_{x^{\ast k}}%
=\delta^{ij}\mathbf{1}.
\end{equation}

\textbf{Creation and annihilation operators}. Let us consider creation and
annihilation operators (CAO) for our case. We introduce creation operators as
the multiplication in $\mathcal{A}$ $a_{x^{i}}^{+}v:=m(x^{i}\otimes v)$ for
every $v\in\mathcal{A}$. We use here the following notation for our CAO
$|i_{1},\ldots,i_{n}\rangle:=m_{n}(x^{i_{1}}\otimes\cdots\otimes x^{i_{n}})$
for state vectors and $a_{x^{i}}\equiv a_{i}^{+}$ for creation operators
corresponding for generators of the algebra $\mathcal{A}$. We have for
example
\begin{equation}%
\begin{array}
[c]{l}%
a_{j_{1}}^{+}\cdots a_{j_{n}}^{+}|0\rangle=|i_{1},\ldots,i_{n}\rangle
.\label{crop}%
\end{array}
\end{equation}
For annihilation operators we assume that $a_{x^{i\ast}}|0\rangle\equiv
a_{i}|0\rangle=0$ for every $x^{\ast i}\in S^{\ast}$ and $a_{s^{\ast}}%
v\in\mathcal{A}^{n-k}$, for $s^{\ast}\in\mathcal{A}^{\ast k}.$The proper
action of annihilation operators on the whole algebra $\mathcal{A}$ is a problem.

If a representation $x$ of quantum logic $L$ in an algebra $A$ is given in
such a way that there is an unique, nondegenerate and positive definite scalar
product then we say that $A$ is the noncommutative Fock space. This means that
our quantum logic are represented by a system with generalized statistics.

\emph{Example 1.} Let $L$ be a quantum logic, where the orthogonality relation
$\perp$ holds for every pair $(\sigma, \sigma^{\prime})$ such that $\ast
\sigma\neq\sigma^{\prime}$. Here the algebra of states $\mathcal{A}$ is the
full tensor algebra $TE$ over the space $E$, and the conjugate algebra
$\mathcal{A}^{\ast}$ is identical with the tensor algebra $TE^{\ast}$. If $T
\equiv0$ then we obtain the most simple example of well--defined system with
generalized statistics. The corresponding statistics is the so--called
infinite (Bolzman) statistics \cite{owg,gre}.

\emph{Example 2} Let $\Psi^{T}$ be a generalized braiding generated by an
operator $T:E^{\ast}\otimes E\rightarrow E\otimes E^{\ast}$. This means that
$\Psi^{T}:TE^{\ast}\otimes TE\rightarrow TE\otimes TE^{\ast}$ is defined as a
set of mappings $\Psi_{k,l}:E^{\ast\otimes k}\otimes E^{\ast\otimes
k}\rightarrow E^{\otimes l}\otimes E^{\ast\otimes k}$, where $\Psi_{1,1}\equiv
R:=T+g_{E}$, and
\begin{equation}%
\begin{array}
[c]{l}%
\Psi_{1,l}:=R_{l}^{(l)}\circ\cdots\circ R_{l}^{(1)},\\
\Psi_{k,l}:=(\Psi_{1,l})^{(1)}\circ\ldots\circ(\Psi_{1,l})^{(k)},\label{up}%
\end{array}
\end{equation}
here $R_{l}^{(i)}:E_{l}^{(i)}\rightarrow E_{l}^{(i+1)}$, $E_{l}^{(i)}%
:=E\otimes\ldots\otimes E^{\ast}\otimes E\otimes\ldots\otimes E$
($l+1$-factors, $E^{\ast}$ on the i-th place, $i\leq l$) is given by the
relation $R_{l}^{(i)}:=\underbrace{id_{E}\otimes\ldots\otimes R\otimes
\ldots\otimes id_{E}}_{l\;\;times}$, where $R$ is on the i-th place,
$(\Psi_{1,l})^{(i)}$ is defined in similar way like $R^{(i)}$. We also
introduce the operator $\tilde{T}:E\otimes E\rightarrow E\otimes E$ by its
matrix elements $(\tilde{T})_{kl}^{ij}=T_{lj}^{ki}$. If the operator
$\tilde{T}$ is a bounded operator acting on some Hilbert space such that we
have the following Yang-Baxter equation on $E\otimes E\otimes E$
\begin{equation}%
\begin{array}
[c]{l}%
(\tilde{T}\otimes id_{E})\circ(id_{E}\otimes\tilde{T})\circ(\tilde{T},\otimes
id_{E})=(id_{E}\otimes\tilde{T})\circ(\tilde{T}\otimes id_{E})\circ
(id_{E}\otimes\tilde{T}),
\end{array}
\end{equation}
and $||\tilde{T}||\leq1$, then according to Bo$\dot{z}$ejko and Speicher
\cite{bs2} there is a positive definite scalar product. Note that the
existence of nontrivial kernel of operator $P_{2}\equiv id_{E\otimes E}%
+\tilde{T}$ is essential for the nondegeneracy of the scalar product
\cite{jswe}. One can see that if this kernel is trivial, then we obtain the
well--defined system with generalized statistics \cite{RM,ral}.

\emph{Example 3:} If a linear and invertible operator $B:E\otimes E\rightarrow
E\otimes E$ defined by its matrix elements $B(x^{i}\otimes x^{j}):=B_{kl}%
^{ij}(x^{k}\otimes x^{l})$ is given such that we have the following
conditions
\begin{equation}%
\begin{array}
[c]{c}%
B^{(1)}B^{(2)}B^{(1)}=B^{(2)}B^{(1)}B^{(2)},\\
B^{(1)}T^{(2)}T^{(1)}=T^{(2)}T^{(1)}B^{(2)},\\
(id_{E\otimes E}+\tilde{T})(id_{E\otimes E}-B)=0,\label{cd}%
\end{array}
\end{equation}
then $I:=gen\{id_{E\otimes E}-B\}$ and one can prove that the corresponding
system is well defined \cite{RM,ral}. In this case
\begin{equation}%
\begin{array}
[c]{l}%
\varrho_{n}^{E}(\tau):=B,
\end{array}
\end{equation}
where $\tau$ is the transposition $(1,2)\rightarrow(2,1)$ in $S_{2}$.

\section{Operads and Fock Spaces \footnotep{Last change May 19, 2003}}

Let $\mathcal{M}(\otimes, \mathbb{C})$ be a monoidal category equipped with a
monoidal operation (a bifunctor) $\otimes: \mathcal{M}\times\mathcal{M}%
\rightarrow\mathcal{M}$, and with the field $\mathbb{C}$ as the unit object.
If there is also a $\ast$-operation $(-)^{\ast}:\mathcal{M}\rightarrow
\mathcal{M}$ and pairing $g = \{g_{U} : U^{\ast}\otimes U\rightarrow
\mathbb{C}\}$ satisfying some known axioms, then we say that we have a
category with duality. One can also include cross or braids, see [Marcinek
1996] for more details.

We would like to construct a new category $\mathcal{P}=\mathcal{P}%
(I,\mathbf{0},\odot,\oplus)$ whose objects are isomorphic classes of of
objects of the initial category $\mathcal{M}$ and equipped with two special
objects $I,\mathbf{0}$, a product $\otimes$ and coproduct $\oplus$ satisfying
the distributivity and associativity conditions
\begin{equation}%
\begin{array}
[c]{ll}%
(\mathcal{U}\odot\mathcal{V})\odot\mathcal{W}\simeq\mathcal{U}\odot
(\mathcal{V}\odot\mathcal{W}) & \mathcal{U}\oplus\mathbf{0}\simeq\mathcal{U}\\
(\mathcal{U}\oplus\mathcal{V})\oplus\mathcal{W}\simeq\mathcal{U}%
\oplus(\mathcal{V})\oplus\mathcal{W}) & \mathcal{U}\odot I\simeq\mathcal{U}\\
\mathcal{U}\odot(\mathcal{V}\oplus\mathcal{W})\simeq(\mathcal{U}%
\odot\mathcal{V})\oplus(\mathcal{U}\odot\mathcal{W}) & \mathcal{U}%
\oplus\mathcal{V}\simeq\mathcal{V}\oplus\mathcal{U}\label{rig}%
\end{array}
\end{equation}
up to natural isomorphisms. We have the following assumptions for our
operadicc description of noncommutative Fock space: 1) $\mathbf{0}$ is the
empty class of objects; 2) $I=\mathbf{1}\mathbb{C}$ is the class corresponding
to the unit object; 3) there is an ordered (finite) collection of classes
$S:=\{x^{i}:i=1,\ldots,N<\infty\}$. These classes form a basis for a finite
linear space $E$ over a field of complex numbers $\mathbb{C}$. There is an
ordered (finite) collection of conjugated classes $S^{\ast}:=\{x^{\ast
i}:i=N,N-1,\ldots,1\}$. They form a basis for the complex conjugate space
$E^{\ast}$. The pairing $(=|-):E^{\ast}\otimes E\longrightarrow Iq$ is given
by $(x^{\ast i}|x^{j}):=\delta^{ij}$. We assume that there is a set of linear
projectors $\odot_{n}:E^{\otimes n}\longrightarrow E^{\otimes n}$, which is
defined uniquely by an iteration procedure
\begin{equation}%
\begin{array}
[c]{c}%
\odot_{n+1}=\odot\circ(id\otimes\odot_{n})=\odot\circ(\odot_{n}\otimes id),
\end{array}
\end{equation}
where $\odot_{1}:=id$, and $\odot_{2}:=\odot$ is given as a data for our
construction. Observe that the definition is unique if and only if our
operation is associative. We are going to describe a collection of spaces
$\mathbb{E}:=\{\mathbb{E}(n):n=0,1,2,\cdots\}$, where $\mathbb{E}(0)=I$,
$\mathbb{E}(1)=E$, and $\mathbb{E}(n):=Im\odot_{n}$ for $n>1$. In this way we
obtain a collection of mappings $\odot_{n}:E^{\otimes n}\rightarrow
\mathbb{E}(n)$, where we use the notation $\odot_{n}(x^{i_{1}}\otimes
\ldots\otimes x^{i_{n}}):=x^{i_{1}}\odot\ldots\odot x^{i_{n}}$. We assume that
the set $\{\odot_{n}(x^{i_{1}},\ldots,x^{i_{n}}):=x^{i_{1}}\odot\ldots\odot
x^{i_{n}}:\sigma=(i_{1},\ldots i_{n})\in\Sigma\}$ forms a basis of
$\mathbb{E}(n)$.

\textbf{Composition operad}. There is a collection of spaces $\mathbb{E}%
:=\{\mathbb{E}(n):n=0,1,2,\cdots\}$, where $\mathbb{E}(0)=I$, $\mathbb{E}%
(1)=E$, and we have the following structure map
\begin{equation}%
\begin{array}
[c]{c}%
\gamma_{\mathbb{E}}:\mathbb{E}(l)\times\mathbb{E}(n_{1})\times\cdots
\times\mathbb{E}(n_{l})\rightarrow\mathbb{E}(n_{1}+\cdots+n_{l})
\end{array}
\end{equation}
for every $n_{1},\ldots,n_{l}=1,2,\ldots$, $l=1,2,\ldots$ such that they form
a nonsymmetric operad. Here we use the following notation $\gamma_{\mathbb{E}%
}(v;v_{1},\cdots,v_{l}):=v(v_{1},\cdots,v_{l})$, where $v\in\mathbb{E}%
(l),v_{1}\in\mathbb{E}(n_{1}),\cdots,v_{l}\in\mathbb{E}(n_{l})$.

\textbf{Coproduct cooperad}. There is the corresponding cooperad
$\mathbb{E}^{\ast}:=\{\mathbb{E}^{\ast}(n):n=0,1,2,\cdots\}$, where
$\mathbb{E}^{\ast}(0)=\mathbb{C}\mathbf{1}^{\ast}$, $\mathbb{E}^{\ast
}(1)=E^{\ast}$, and with the following structure map
\begin{equation}%
\begin{array}
[c]{c}%
\gamma_{\mathbb{E}^{\ast}}:\mathbb{E}^{\ast}(l)\times\mathbb{E}^{\ast}%
(n_{1})\times\cdots\times\mathbb{E}^{\ast}(n_{l})\rightarrow\mathbb{E}^{\ast
}(n_{1}+\cdots+n_{l})
\end{array}
\end{equation}
for every $n_{1},\ldots,n_{l}=1,2,\ldots$, $l=1,2,\ldots$. We define creation
operators $a_{x^{i}}^{+}v:=\odot_{n+1}(x^{i}\otimes v)$, where $x^{i}\in E$,
$v\in\mathbb{E}(n)$, and annihilation ones $a_{x^{\ast i}}^{-}\,x^{j_{1}}%
\odot\ldots\odot x^{j_{n}}:=(x^{\ast i}\mid x^{j_{1}})_{1}\;x^{j_{2}}%
\odot\ldots\odot x^{j_{n}}$. We have here the following relations
\begin{equation}%
\begin{array}
[c]{c}%
a_{\ast i}^{-}\circ a_{j}^{+}=(x^{\ast i}\mid x^{j})_{1}.
\end{array}
\end{equation}
For every $x^{\ast i}$ we define an operator $\iota_{x^{\ast i}}$ by the
formula $\iota_{x^{\ast i}}x^{j}:=(x^{\ast i}\mid x^{j})_{1}$. We use the
notation
\begin{equation}
\iota_{x^{\ast i}}^{(k)}x^{j_{1}}\odot\ldots\odot x^{j_{k}}\odot\ldots\odot
x^{j_{n}}=x^{j_{1}}\odot\ldots\odot\underbrace{\iota_{x^{\ast i}}}_{k}%
x^{j_{k}}\odot\ldots\odot x^{j_{n}}.
\end{equation}
We introduce new operators
\begin{equation}%
\begin{array}
[c]{c}%
b_{i}^{+}:=a_{i}^{+},\\
b_{\ast i}^{-}(x^{j_{1}}\odot\ldots\odot x^{j_{k}}\odot\ldots\odot x^{i_{n}%
})=\sum_{k=1}^{n}\left(  \iota_{x^{\ast i}}^{(k)}\circ\mathbb{T}(k)\right)
x^{\ast i}\odot x^{j_{1}}\odot\ldots\odot x^{j_{k}}\odot\ldots\odot x^{j_{n}}.
\end{array}
\end{equation}
We obtain here the following relations
\begin{equation}%
\begin{array}
[c]{c}%
b_{\ast i}^{-}\circ b_{i}^{+}-\sum_{k,l}T_{kl}^{ij}\;b_{k}^{+}b_{\ast l}%
^{-}=(x^{\ast i}\mid x^{j}),
\end{array}
\end{equation}
where $T_{kl}^{ij}$ are matrix elements $\mathbb{T}(1)(x^{\ast i}\otimes
x^{j}):=\sum_{k,l}T_{kl}^{ij}x^{k}\otimes x^{\ast l}$. Paring operad is a
collection of pairings $(-|-):=\{(-|-)_{n}:n=1,2,\ldots\}$, where
$(-|-)_{0}:=0$, $(-|-)_{1}$ is given above. Wick ordering operad is a
collection of operators $\{\mathbb{T}(n):\mathbb{E}^{\ast}(1)\otimes
\mathbb{E}(n)\rightarrow\mathbb{E}(n)\otimes\mathbb{E}^{\ast}(1)\}$.

\textbf{Fock operad representation}. There are two familes of operators
\begin{equation}%
\begin{array}
[c]{c}%
a^{+} := \{a^{+}_{i}:\mathbb{E}(n)\rightarrow\mathbb{E}(n+1), i=1,\ldots, N\},
\end{array}
\end{equation}
and
\begin{equation}%
\begin{array}
[c]{c}%
a^{-} := \{a^{-}_{\ast i}:\mathbb{E}(n)\rightarrow\mathbb{E}(n-1), \ast i=N,
N-1,\ldots, 1 \}.
\end{array}
\end{equation}

\textbf{Noncommutative Fock space}. If there is a set of nondegenerate and
positive definite scalar products $\langle-|-\rangle:= \{\langle-|-\rangle_{n}
: n=1, 2, \ldots\}$ such that $\langle0|0 \rangle= 0$, $\langle x^{i}%
|x^{j}\rangle_{1} := (x^{\ast i}|x^{j})$, $\langle u|v\rangle_{n} := (u^{\ast
}|v)_{n}$ and
\begin{equation}%
\begin{array}
[c]{cc}%
a^{-}_{x^{\ast i}} |0\rangle= 0, & \langle a^{+}_{i} v|w \rangle= \langle
v|a^{-}_{\ast i}w\rangle,
\end{array}
\end{equation}
then we say that our quantum system is well defined.

\textbf{Commutation relations}. Let us describe a series of representations
\begin{equation}%
\begin{array}
[c]{c}%
\varrho_{n}:\pi\in S_{n}\mapsto S_{\pi}\in\mathrm{End}(E^{\otimes n}),
\end{array}
\end{equation}
where $S_{n}$ is the symmetric group, and $S_{\pi}:=S_{\tau_{i_{1}}}%
\circ\ldots S_{\tau_{i_{k}}}$, corresponds for the following decomposition of
$\pi$ on transpositions $\pi:=\tau_{i_{1}}\circ\ldots\circ\tau_{i_{k}}$. If
$\odot_{n}$ is a composition mapping and $\varrho_{n}$ is a representation of
the symmetric group $S_{n}$ described above, then for every permutation
$\pi\in S_{n}$ there is a new composition mapping $\odot_{n}\circ S_{\pi}$.

\section{Particle processes and cobordisms with trees
\footnotep{Last change March 1, 2000}}

Let $\mathcal{C} := \{\Sigma\}$ be a collection of $d$--dimensional compact,
oriented and smooth manifolds without boundary. All these manifolds can be
multiconnected in general. If $\Sigma_{1} , \Sigma_{2}\in\mathcal{C}$, then we
denote by $\Sigma_{1}\cup\Sigma_{2}$ their disjoint sum. For every manifold
$\Sigma\in\mathcal{C}$ there is the corresponding manifold $\Sigma^{\ast}$
such that $\partial(\Sigma\times[0,1]) = \Sigma\cup\Sigma^{\ast}$. We also
assume for our study that two finite set $\Gamma:= \{\mathbf{i}_{1} , \cdots,
\mathbf{i}_{N}\}$ and $\Gamma^{\ast} := \{\mathbf{i}^{\ast}_{N} , \cdots,
\mathbf{i}^{\ast}_{1}\}$ are given.

Let us assume that there is a corresponding collection $\mathcal{C}%
(\mathcal{Q}):=\{(\Sigma,\mathcal{Q}):\Sigma\in\mathcal{C},\mathcal{Q}%
\in\mathcal{Q}(\Gamma)\}$ of manifolds with discrete structure $\mathcal{Q}%
(\Gamma):=\{\mathcal{Q}(n):n=0,1,\ldots\}$ on these manifolds, where
$\mathcal{Q}(n)$ is a finite set of points of $\Sigma$ with labels in $\Gamma
$, $\mathcal{Q}(0)=\emptyset$. The corresponding collection $\mathcal{C}%
(\mathcal{Q}^{\ast}):=\{(\Sigma^{\ast},\mathcal{Q}^{\ast}):\Sigma^{\ast}%
\in\mathcal{C},\mathcal{Q}^{\ast}\in\mathcal{Q}(\Gamma^{\ast})\}$ is defined
in a similar way. In our approach particle processes of interactions are
described by pairs of the form $(\mathcal{M},G)$ which transforms the initial
configuration $(\Sigma_{0},\mathcal{Q}_{0})$ into the outgoing one
$(\Sigma_{1},\mathcal{Q}_{1})$ representing the results of interactions. Here
$\mathcal{M}:\Sigma_{0}\longrightarrow\Sigma_{1}$ transforms the manifold
$\Sigma_{0}$ into $\Sigma_{1}$, and $G:\mathcal{Q}_{0}\longrightarrow
\mathcal{Q}_{1}$ transforms the corresponding discrete structures. If
$(\mathcal{N},G^{\prime})$ is the second pair, then we assume that there is a
pair $(\mathcal{M}\circ\mathcal{N},G\circ G^{\prime})$, the composition of
$(\mathcal{M},G)$ and $(\mathcal{N},G^{\prime})$. We can use the concept of
cobordisms manifold and rooted trees for the description of these mappings and
their compositions. We describe here an arbitrary particle process as a
cobordism manifold with a tree structure. Denote by $\mathcal{E}$ a collection
of $d+1$--dimensional compact, oriented and smooth manifolds with boundary. We
assume that for every manifold $\mathcal{M}\in\mathcal{E}$ with boundary
$\partial\mathcal{M}$ there are two manifolds $\Sigma_{0}$ and $\Sigma_{1}$ in
$\mathcal{C}$ such that the boundary $\partial\mathcal{M}$ is diffeomorphic to
$\Sigma_{0}\cup\Sigma_{1}^{\ast}$.

A $d+1$ manifold $\mathcal{M}\in\mathcal{E}$ with boundary $\partial
\mathcal{M}$ such that there are two smooth diffeomorphisms $f_{0}:\Sigma
_{0}^{\ast}\longrightarrow\partial\mathcal{M}$, and $f_{1}:\Sigma
_{1}\longrightarrow\partial\mathcal{M}$ is said to be a cobordism of
$\Sigma_{0}$ and $\Sigma_{1}$ and it is denoted by $\mathcal{M}(f_{0},f_{1})$.
Two cobordisms $\mathcal{M}(f_{0},f_{1})$ and $\mathcal{M}^{\prime}%
(f_{0}^{\prime},f_{1}^{\prime})$ are said to be equivalent if there is a
diffeomorphism $F:\mathcal{M}\longrightarrow\mathcal{M}^{\prime}$ such that
$f_{0}^{\prime}=Ff_{0}$ and $f_{1}^{\prime}=Ff_{1}$. The equivalence class of
cobordisms of $\Sigma_{0}$ and $\Sigma_{1}$ up to diffeomorphisms is denoted
by $\mathcal{M}(\Sigma_{0},\Sigma_{1})$. The collection of all classes
cobordisms of $\Sigma_{0}$ and $\Sigma_{1}$ is denoted by $\mathcal{E}%
(\Sigma_{0},\Sigma_{1})$. Let $\mathcal{N}$ be a cobordism of $\Sigma_{2}$ and
$\Sigma_{3}$ with diffeomorphisms $f_{0}^{\prime}:\Sigma_{2}^{\ast
}\longrightarrow\partial\mathcal{N}$ and $f_{1}^{\prime}:\Sigma_{3}%
\longrightarrow\partial\mathcal{N}$. In certain cases we can glue two
cobordisms $\mathcal{M}$ and $\mathcal{N}$ along $\Sigma_{1}$ and $\Sigma_{2}$
by identifying the part of boundary of $\mathcal{M}$ diffeomorphic to
$\Sigma_{1}$ with the part of $\partial\mathcal{N}$ diffeomrphic to
$\Sigma_{2}^{\ast}$, respectively. The composition $(f_{0}^{\prime})^{-1}%
f_{1}$ is then a diffeomorphism of $\Sigma_{1}$ onto $\Sigma_{2}^{\ast}$. In
this way we obtain the cobordism of $\Sigma_{0}$ and $\Sigma_{2}$. We denote
it by $\mathcal{M}\circ\mathcal{N}$. The operation of gluing of cobordisms up
to diffeomorphisms define the following composition $\circ:\mathcal{E}%
(\Sigma_{0},\Sigma_{1})\times\mathcal{E}(\Sigma_{1},\Sigma_{2})\longrightarrow
\mathcal{E}(\Sigma_{0},\Sigma_{2})$. One can glue three or more cobordimsms.
One can see that these gluings define a semigroups structure on the collection
of all cobordism for certain class of manifolds.

It is interesting to study cobordisms of manifolds with certain additional
structures. If for example these manifolds are configuration spaces for
multiparticle system, then we must restrict our gluings for those which are
admissible for our physical problem. For this goal we need some additional
assumptins. A \textit{rooted tree} is a finite, loop free, connected graph
which contains edges and nodes. There is one distinguisched edge called a
root. Every node is inner and there are no outer nodes. A rooted tree with one
node and $n$ entrance edgees is said to be a \textit{the prime n--tree}. In
our physical interpretation entrance edges describe initial configurations of
particles. The root is said to be an exit. It represents the unique final
configuration. There also the corresponding concept of co--rooted trees.
Starting with a set of prime $2$--trees and $2$--co--trees one can construct a
graph $G$ with arbitrary number of entrance and outgoing edges. Such graph can
represent arbitrary particle processes. If we embed trees into the cobordism
manifold, then we obtain cobordisms with a tree structure. More precisely, let
$\mathcal{M}(\Sigma_{0},\Sigma_{1})$ be a cobordism of $\Sigma_{0}$ and
$\Sigma_{1}$ with discrete structures $\mathcal{Q}_{0}:=\mathcal{Q}(n)$ and
$\mathcal{Q}_{1}:=\mathcal{Q}(m)$, respevtively, and let $G$ be a graph with
$n$ entrance edges and $m$ outgoing ones, then for the embeding of $G$ in
$\mathcal{M}$ we must connect $n$ entrance edges of $G$ with $n$ points of
$\mathcal{Q}_{n}(\Sigma_{0})$, and $m$ outgoing edges with $m$ points of
$\mathcal{Q}_{1}$. Note that we need here some additional restrictions for the
equivalence of cobordisms. If the disrete structure represents fermions and
bosons, then we need certain supermanifolds or some generalization for to
describe the corresponding cobordisms. We can use here the notion of the
so--called semisupermanifolds introduced by S. Duplij \cite{sdu}.

We denote by $Cob=Cob_{d}(\mathcal{Q})$ the category of cobordisms with a tree
structure. Objects of this category are $d$--dimensional compact, oriented and
smooth manifolds $\Sigma$ equipped with a discrete structure $\mathcal{Q}%
(\Sigma)$. Morphisms are cobordism manifold with a tree structure. Composition
of morphisms can be expressed as a gluings of cobordisms with trees.

Let us assume that we have the so--called multiparticle states operad
$\Lambda:=\{\Lambda(n):n=0,1,2,\cdots\}$, where $\Lambda(n)$ is an
$n$--particle space of states, $\Lambda(0)=\mathbb{C}\mathbf{1}$, $\Lambda(1)$
is equipped with a basis with $N$ elements $x^{1},\cdots,x^{N}$, $\Lambda(n)$
for $n\geq2$ can be obtained from $\Lambda_{1}$ by an $n$-ary operation
$\circ_{n}:\Lambda\times\cdots\times\Lambda\longrightarrow\Lambda(n)$. If the
$n$-ary operation can be obtained uniquely from a binary one $F\equiv
F_{2}:\Lambda\times\Lambda\longrightarrow\Lambda(2)$ by an iteration
procedure, then we say that the configuration operad is well defined. One can
introduce a category $Part\equiv Part(\Lambda)$ related to the operad
$\Lambda$. Let us consider an arbitrary functor $\mathcal{Z}%
:Cob\longrightarrow Part$. Let $\mathcal{M}$ be a cobordism of $\Sigma_{0}$
and $\Sigma_{1}$, then the our goal is the construction of $\mathcal{Z}%
(\Sigma_{0})$, $\mathcal{Z}(\Sigma_{1})$, and $\mathcal{Z}(\mathcal{M})$ as an
mapping $\Phi_{\mathcal{M}}:\mathcal{Z}(\Sigma_{0})\longrightarrow
\mathcal{Z}(\Sigma_{1})$. We express $\Phi$ as the path integral
\begin{equation}
(\Phi_{\mathcal{M}}\mathcal{W})(\varphi)=\int\exp(-S(\varphi))\mathcal{W}%
(\varphi)D\varphi,
\end{equation}
where $\varphi$ is a field on $\mathcal{M}$ with a given boundary conditions
on $\Sigma_{0}$ and $\Sigma_{1}$, $S$ is a given action and $\mathcal{W}$ is
an observable. Note that $\mathcal{Z}(\Sigma_{k})$ should be expressed as a
sequence $\Lambda(n_{1})\times\cdots\times\Lambda(n_{r})$, where $\Sigma
_{i=1}^{r}n_{i}=n$ is the number of point particles represented by
$\mathcal{Q}(\Sigma_{k})$.

 \end{document}